\documentclass{mcom-l}

\newcommand{\Z}{\mathbb{Z}}
\newcommand{\N}{\mathbb{N}}

\newcommand{\lp}{\mathcal{F}\backslash\mathcal{S}}
\newcommand{\norm}[1]{\left\lVert#1\right\rVert}

\usepackage{amssymb}
\usepackage{graphicx}
\usepackage{booktabs}
\usepackage{tabularx}
\usepackage{adjustbox}
\usepackage{hyperref}
\usepackage{cleveref}
\usepackage{algorithm}
\usepackage{algpseudocode}
\usepackage{float}
\usepackage{caption}
\usepackage{subcaption}
\usepackage{enumitem}   
\usepackage{alltt}
\algtext*{EndWhile}
\algtext*{EndFor}
\algtext*{EndIf}

\copyrightinfo{}{}

\newtheorem{theorem}{Theorem}[section]

\theoremstyle{definition}

\theoremstyle{remark}

\numberwithin{equation}{section}

\begin{document}

\title[Smooth Subsum Search]{Smooth Subsum Search\\ A heuristic for practical integer factorization}


\author[M. Hittmeir]{Markus Hittmeir}
\address{SBA Research, Floragasse 7, 1040 Vienna, Austria}
\curraddr{}
\email{mhittmeir@sba-research.org}
\thanks{SBA Research (SBA-K1) is a COMET Centre within the framework of COMET – Competence Centers for Excellent Technologies Programme and funded by BMK, BMDW, and the federal state of Vienna. The COMET Programme is managed by FFG}


\date{}

\dedicatory{}

\begin{abstract}
The two currently fastest general-purpose integer factorization algorithms are the Quadratic Sieve and the Number Field Sieve. Both techniques are used to find so-called smooth values of certain polynomials, i.e., values that factor completely over a set of small primes (the factor base). As the names of the methods suggest, a sieving procedure is used for the task of quickly identifying smooth values among the candidates in a certain range. While the Number Field Sieve is asymptotically faster, the Quadratic Sieve is still considered the most efficient factorization technique for numbers up to around 100 digits. In this paper, we challenge the Quadratic Sieve by presenting a novel approach based on representing smoothness candidates as sums that are always divisible by several of the primes in the factor base. The resulting values are generally smaller than those considered in the Quadratic Sieve, increasing the likelihood of them being smooth. Using the fastest implementations of the Self-initializing Quadratic Sieve in Python as benchmarks, a Python implementation of our approach runs consistently 5 to 7 times faster for numbers with 45-100 digits, and around 10 times faster for numbers with 30-40 digits. We discuss several avenues for further improvements and applications of the technique.
\end{abstract}

\maketitle

\section{Introduction}
Integer factorization is the task of computing the divisors of natural numbers. It is a problem with a long and fascinating history, and it is certainly among the most influential in algorithmic number theory. While there is a variety of algorithms significantly faster than the brute-force search for divisors, it is still an open problem to construct a technique that efficiently factors general numbers with hundreds to thousands of digits. The hardness of this problem is fundamental for the security of widely used cryptographical schemes, most prominently the RSA cryptosystem. Nevertheless, there is no proof for its hardness besides the fact that decades of efforts have failed to construct a more efficient technique.

Quite regularly, there are set new records\footnote{A recent record was the factorization of a 250 decimal digit number in February 2020 by Boudot et al., see \url{https://en.wikipedia.org/wiki/RSA_Factoring_Challenge}} concerning the factorization of numbers of certain size, mostly due to improved implementations of the best available algorithms and advances in the hardware and computing power. In addition, the bound for the deterministic integer factorization problem has been improved multiple times in recent years (\cite{Hit0}, \cite{Hit1}, \cite{Har}, \cite{HH}). On the other hand, there has only been little progress in the development of new techniques for practical integer factorization since the invention of the Number Field Sieve (\cite{Len1}) in the 1990s. One of the earlier algorithms with sub-exponential runtime was by Dixon (\cite{dix}) in 1981. From today's perspective, it may be considered as a prototype for several other algorithms. To this group belong the Continued Fraction factorization method (CFRAC) described and implemented by Morrison and Brillhart (\cite{mor}) in 1975, the Linear Sieve by Schroeppel and the Quadratic Sieve by Pomerance. The latter author analyzed and compared these algorithms in 1982 in \cite{pom}. In general, there is an extensive amount of literature on practical integer factorization algorithms. The reader may find information on the mentioned methods and other factorization techniques in the survey \cite{Len} and in the monographs \cite{Rie} and \cite{Wag}. These sources also discuss a variety of factorization algorithms for numbers $N$ that satisfy certain properties, or that have prime factors of certain shape. For example, Fermat's factorization method runs particularly fast if $N$ has co-divisors that are very close, while Pollard's $P-1$ technique and the Elliptic Curve Method (ECM) work particularly well for numbers with small prime factors. In the present paper, we are interested in so-called general-purpose factorization algorithms, which means that the runtime complexity estimate of the procedure only depends on the size of the input number $N$. Both the Quadratic Sieve and the Number Field Sieve belong to this group. From a practical point of view, the Quadratic Sieve (and its modifications) is considered the best general-purpose factorization algorithm for $N$ up to around 100 digits, whereas the Number Field Sieve is faster for inputs beyond that. While methods such as CFRAC and ECM stay somewhat competitive up to a small input length (see \cite{mil}), factorization tools usually use ECM to rule out the existence of relatively small prime factors, but then switch either to the Number Field Sieve or to the Self-initializing Quadratic Sieve (SIQS), which is the fastest modification of the original Quadratic Sieve algorithm.\footnote{In its standard settings, the tool \emph{Msieve} (\url{https://sourceforge.net/projects/msieve/}) uses ECM for searching for small prime factors with up to 15 digits. Similarly, for numbers with 31 digits or more, the calculator of Dario Alpern (\url{https://www.alpertron.com.ar}) switches to SIQS after a quick search for small factors via ECM.}

Our contribution is the presentation and evaluation of a novel approach, \emph{Smooth Subsum Search} (SSS). Just like CFRAC, the Linear Sieve and the Quadratic Sieve, SSS belongs to the group of Dixon-type algorithms. The core idea is based on constructing smoothness candidates in a way that guarantees divisibility by many primes, while keeping their values as small as possible. This is achieved by finding collisions of affine transformations of the roots of the considered polynomial modulo the primes in the factor base. The construction allows to quickly generate values that have a high likelihood of being smooth without sieving the complete interval for each prime of the factor base. In addition to a detailed explanation of our algorithm, a discussion of its advantages compared to the Quadratic Sieve and theoretical arguments for the runtime improvement, we will consider our current implementation in Python and several experiments on the runtime complexity of SSS. Our results demonstrate that SSS outperforms SIQS on semiprime\footnote{Throughout the paper, we will use this term for numbers composed of two distinct prime factors. Under certain conditions, such numbers are the most difficult to factorize.} numbers with 30 to 100 decimal digits, finding smooth numbers around 5 to 7 times faster. While these findings need to be corroborated in experiments with different implementations of SIQS in other programming languages, they indicate that SSS is the fastest available factorization algorithm for general numbers in this range. 

The remainder of the paper is structured as follows: In \Cref{sec:rel_work}, we discuss the Quadratic Sieve and its modifications in greater detail, together with an important subroutine of our own algorithm. \Cref{sec:sss} presents the full algorithm SSS as it is currently implemented. The results of our experiments can be found in \Cref{sec:experiments}, where SSS is compared to the fastest implementations of the Quadratic Sieve in Python. \Cref{sec:improvements} explores three speculative ideas on how to construct further variants and improvements of SSS. Finally, we summarize and give concluding remarks in \Cref{sec:conclusion}.

\section{Related Work} \label{sec:rel_work}
Let $N$ be the number we want to factorize. We will always assume that $N$ is odd, composite and not a perfect power of another number. The currently fastest general-purpose factorization algorithm is the General Number Field Sieve with a heuristic asymptotic runtime complexity of 
\[
\exp\left((\sqrt[3]{64/9}+o(1))(\log N)^{1/3}(\log\log N)^{2/3}\right).
\]
For comparison, the heuristic asymptotic runtime complexity of the Quadratic Sieve is $\exp\left((1+o(1))(\log N)^{1/2}(\log\log N)^{1/2}\right)$ (\cite{pom1}). The Number Field Sieve is significantly more complicated and harder to implement than the Quadratic Sieve. As mentioned, the latter technique is also still faster in practice for factoring numbers up to around 100 digits. We now restrict our attention to the technical details of the Quadratic Sieve, as it is a suitable representative for the other Dixon-type methods. In addition, it will be used as a benchmark for our own algorithm in the experiments in \Cref{sec:experiments}. 

The ultimate goal of the Quadratic Sieve and the other Dixon-type algorithms is to find two numbers $X$ and $Y$ such that $X^2\equiv Y^2\mod N$. For such pairs, one can easily show that there is a high chance that $\gcd(X-Y,N)$ gives a nontrivial divisor. However, perfect squares are a rare occurrence. Since there are roughly $\sqrt{N}$ square numbers below $N$, simply picking random values for $X$ and hoping for $X^2\pmod N$ to be square will take too long. So in order to find such $X$ and $Y$, the first phase of all these algorithms tries to collect so-called \emph{smooth relations}. Such a relation is a congruence $x^2\equiv y \mod N$ such that the prime factorization of $y$ only consists of primes up to a certain bound $B$. In this case, the number $y$ is called \emph{$B$-smooth} (or just $\emph{smooth}$). Let $p_1,\ldots, p_m$ be the primes up to $B$. The set $\{p_1,\ldots, p_m\}$ of these primes is referred to as \emph{factor base}. The smooth relations may then be written as
$
x^2 \equiv y = p_1^{e_{1}}\cdots p_m^{e_{m}}\mod N
$
for certain exponents $e_i\in\N_0$, $i=1,\ldots,m$.
The first phase proceeds until enough of these relations are found, which usually means a bit more than the number $m$ of primes in the factor base. As soon as this is the case, the second phase starts. Here, we want to find linear dependencies in the exponent vectors $(e_1,\ldots, e_m)$ of the smooth relations, i.e., a subset of these vectors that sum to the $0$-vector. We can use well-known methods from linear algebra to identify such subsets and thereby find smooth relations that can be multiplied to obtain a new congruence of the form 
\[
x_1^2x_2^2\cdots x_k^2\equiv \prod_{i=1}^m p_i^{e_{i,1}+e_{i,2}+\cdots+e_{i,k}} \mod N,
\]
where $e_{i,1}+e_{i,2}+\cdots+e_{i,k}\equiv 0 \mod 2$ for all $i$. As a consequence, there exists a number $Y$ such that $Y^2=\prod_{i=1}^m p_i^{e_{i,1}+e_{i,2}+\cdots+e_{i,k}}$. Clearly, there also exists $X$ such that $X^2=x_1^2x_2^2\cdots x_k^2$, and we have found what we were looking for.
One important parameter of these factorization techniques is the smoothness bound $B$. If $B$ is very large, then it is easier to find smooth relations. On the other hand, we need to find much more of them for the second phase to work. In addition, the matrix in the linear algebra step will be very large. If $B$ is rather small, then we will most likely not be able to find enough smooth relations.

The main difference between Dixon's method, CFRAC, the Quadratic Sieve and our algorithm SSS is in how the smooth relations are collected in the first phase. Dixon's algorithm chooses random values $x$ and checks if $x^2\pmod N$ is $B$-smooth. CFRAC uses the convergents in the continued fraction expansion of $\sqrt{N}$. While the value of $x^2\pmod N$ in Dixon's method is bounded by $N$, the main advantage of CFRAC over Dixon is that it generates congruences of the form $x^2\equiv y \mod N$ with $|y|\leq 2\sqrt{N}$, which increases the chances for $y$ being smooth. In the basic version of the Quadratic Sieve, the values of the polynomial $f(x)=(x+\lceil \sqrt{N} \rceil)^2-N$ are used as candidates for possibly smooth $y$. Note that 
\begin{equation}\label{eq:bound}
y=f(x)=x^2 + 2x\lceil \sqrt{N}\rceil + (\lceil \sqrt{N}\rceil)^2-N,
\end{equation}
which for small $x$ gives a bound similar to the one obtained in CFRAC. The main advantage of the Quadratic Sieve compared to CFRAC is in how we check the smoothness of the candidates $y$. While CFRAC applies trial division or techniques for smoothness detection of batches of numbers (see \cite{Ber}), the Quadratic Sieve uses \emph{sieving}. The basic idea is as follows: Let again $\{p_1,\ldots,p_m\}$ be the factor base $\mathcal{F}$. Our goal is to find arguments $x$ such that $f(x)$ is divisible by many primes in $\mathcal{F}$. We first remove all primes $p$ from $\mathcal{F}$ for which Legendre's symbol $(N|p)=-1$, since $N$ is a quadratic non-residue modulo such $p$ and, thus, $p$ cannot divide $f(x)$ for any value of $x$. Next, we use the well-known Tonelli-Shanks algorithm to compute the solutions of $f(x)\equiv 0 \mod p$ for each remaining prime $p$ in $\mathcal{F}$. Note that 
\[
f(x+kp)\equiv f(x) \mod p 
\]
for every $k\in\Z$. Hence, if $f(x_0)$ is divisible by $p$, so is $f(x_0+kp)$ for every $k\in\Z$. 
We start the sieve by initializing an indexed array of pre-defined length $L$ containing zeros, where each index corresponds to an argument $x$. For each $p\in\mathcal{F}$ and each $x_0$ with $f(x_0)\equiv 0 \mod p$, we then raise the entries in the array at the indices $x_0+kp$ for $k\in\Z$ by $\lceil\log p\rceil$. At the end of the sieving, we check the array for entries that are larger than a certain threshold, as for these it follows that $f(x)$ is divisible by a lot of primes in the factor base and thus likely to be smooth. Since the information about the actual divisors of such $f(x)$ is lost, we still need to apply either trial division or some other factorization algorithm suitable for finding small prime factors. However, the sieving procedure greatly reduces the total number of such applications compared to what is needed in CFRAC. 

If we have not found enough smooth relations by applying the sieving procedure above, we may continue by increasing the length $L$ of the sieving array. However, this will lead to the expression (\ref{eq:bound}) getting quite large, which decreases the chances of $f(x)$ being smooth. The Multiple Polynomial Quadratic Sieve (MPQS) (\cite{Sil}), an improvement of the basic version of the quadratic sieve, makes use of multiple polynomials $f_{a,b}(x)=(ax+b)^2-N$ with $a,b\in\Z$. While switching to another polynomial allows the use of a new array and thus solves the problem of increasing values of the smoothness candidates, it is also somewhat expensive to initialize new polynomials and compute the solutions of $f_{a,b}(x)\equiv 0\mod p$ for all the primes $p$ in the factor base. A remedy for this problem has been found in another improvement, the already mentioned Self-initializing Quadratic Sieve (SIQS) \cite{Con}. The idea is to choose $a$ as a product of primes in the factor base, and $b$ such that $b^2-N=ac$ for some $c\in\Z$. Then $f_{a,b}(x)=a(ax^2+bx+c)$, so $f_{a,b}(x)$ is smooth if and only if $ax^2+bx+c$ is smooth. In addition, there are several different possible choices $b$ for each choice of $a$ that may be computed via the Chinese Remainder Theorem. The corresponding solutions to $f_{a,b}(x)\equiv 0\mod p$ for each prime $p$ are related in a way that allows to change to a new polynomial much more easily. In particular, the expensive inversion of elements modulo the primes $p$ only needs to be done once for every choice of $a$.

One further improvement that may be applied to all Dixon-type algorithms is the so-called \emph{large prime variant} (\cite{Boe}). If we have found a relation $x^2\equiv y\mod N$ where $y$ is partially smooth, meaning that the divisor $r$ in its factorization that is not smooth is comparably small, then we may also save this so-called \emph{partial relation}. As soon as we find two partial relations of the form
\begin{align*}
x_1^2\equiv y_1\cdot r \mod N \\
x_2^2\equiv y_2\cdot r \mod N 
\end{align*}
where $y_1$ and $y_2$ are smooth, we may combine them to get a full relation, namely $(r^{-1})^2x_1^2x_2^2\equiv y_1y_2\mod N$. This improvement can speed up the search for smooth relations by about a factor of $2$.

Finally, let us briefly discuss the smoothness detection algorithm by Bernstein that has already been mentioned in the context of CFRAC. For cases where sieving the values of a polynomial in the described sense is not an option, Bernstein describes another way to find smooth and partially smooth numbers in a given set. As mentioned in the abstract, SSS is not based on sieving. In \Cref{sec:sss}, we will use the following procedure (\cite[Algorithm 2.1]{Ber}) instead. Let $\mathcal{F}=\{p_1,\ldots,p_m\}$ be the factor base and $\mathcal{C}=\{x_1,\ldots,x_n\}$ be a set of positive integers.
\begin{algorithm}\caption*{\textbf{smooth\_batch}($\mathcal{F}$, $\mathcal{C}$)}
\begin{algorithmic}[1]
		\State{Compute $z=p_1\cdots p_m$ using a product tree.}
		\State{Compute $z\pmod{x_1},\ldots, z\pmod {x_m}$ using a remainder tree.}
		\State{For $k=1,\ldots,n$, compute $y_k=(z\pmod{x_k})^{2^e}\pmod{x_k}$ using repeated squaring. Here, $e$ is the smallest nonnegative integer such that $2^{2^e}\geq x_k$.}
		\State{For $k=1,\ldots,n$, print $x_k/\gcd(x_k,y_k)$.}
\end{algorithmic}
\end{algorithm}

The algorithm outputs the part of the factorization of each $x_k$ that is not smooth.
For more information on smooth detection as well as on product and remainder trees, see \cite[Section 3.3]{Cra}. We now move on to the presentation of SSS.

\section{Smooth Subsum Search}\label{sec:sss}

Let us consider an integer polynomial $f:\Z \rightarrow \Z$ and a factor base $\mathcal{F}$ of primes smaller than some bound $B$. The task is to find $x\in\Z$ such that $f(x)$ is $B$-smooth. In \Cref{sec:rel_work}, we have discussed the Quadratic Sieve for solving this problem for the polynomial $f(x)=(x+\lceil \sqrt{N} \rceil)^2-N$ in order to factor the number $N$. The algorithm presented in this section solves the same problem on the same polynomial, but uses a heuristic search approach instead of sieving. 

\subsection{The core idea}\label{sec:bi}
We start by discussing a suitable representation of the possible solutions $x$. We will see that this representation already restricts our attention to values $x$ with increased likelihood of $f(x)$ being smooth.
We are looking for preferably small values $x$ such that $f(x)$ is divisible by several (powers of) the primes in $\mathcal{F}$. Let $p_1,p_2,\ldots,p_m$ be the odd primes in $\mathcal{F}$, and assume that $(N|p_i)=1$ for $i=1,\ldots,m$.  Then, the congruence 
$
f(x)\equiv 0 \mod p_i
$
has two solutions $s_{i,1}, s_{i,2}$. Using the Chinese Remainder Theorem, we may easily compute those $x$ for which $f(x)$ is divisible by specific primes in $\mathcal{F}$ in form of sums (see \textbf{get\_x} below). Namely, with $x_i\in \{0,1,2\}$, we represent the candidates $x$ via vectors 
\[
x\,\widehat{=}\, (x_1,x_2,\ldots,x_m).
\]
If $x_i>0$, we fix the value of $x\pmod{p_i}=s_{i,x_i}$ and thereby make sure that $f(x)\equiv 0 \mod p_i$ holds. If $x_i=0$, we leave the value of $x\pmod{p_i}$ undefined, and it will ultimately be determined by the other values in the vector that we did fix.
Of course we would like $f(x)$ to be divisible by as many primes in $\mathcal{F}$ as possible. However, if we run the Chinese Remainder Theorem on vectors $(x_1,x_2,\ldots,x_m)$ where most of the $x_i$ are nonzero, then the modulus $M$, i.e. the product of the primes $p_i$ for which $x\pmod{p_i}$ is fixed, will be very large, and so will be the average values of the resulting sums $x$. The chances for $f(x)$ being smooth are then rather small. From this observation, it becomes quite clear that we would like to keep $|x|/M$ as small as possible in order to increase our chances for finding smooth numbers. What we are looking for are sparse solution vectors leading to small \emph{subsums} $x$. We therefore restrict our attention to a subset $\mathcal{S}$ of the smaller primes in $\mathcal{F}$ and, for generating a random solution, select only a few primes in $\mathcal{S}$ for which the corresponding value $x_i$ is nonzero. We call $\mathcal{S}$ the \textit{small factor base}. Assume that there are $n$ primes in $\mathcal{S}$, i.e., let $\mathcal{S}=\{p_1,\ldots,p_n\}$. Then $x_i=0$ for $i>n$ in all considered representations $(x_1,x_2,\ldots,x_m)$, and we may just denote them by $(x_1,x_2,\ldots,x_n)$. 
Here is our function for retrieving the value $x$ that corresponds to the representations $(x_1,x_2,\ldots,x_n)$ by applying the Chinese Remainder Theorem.
\begin{algorithm}
	\caption*{\textbf{get\_x$(x_1,x_2,\ldots,x_n)$}}
	\begin{algorithmic}[1]
		\State{Compute $M=\prod_{i:\, x_i\neq 0}p_i$.}
		\State{Return 
			\[
			x:=\sum_{i:\, x_i\neq 0} \frac{M c_i}{p_i}\cdot s_{i,x_i} \pmod M,
			\]
		where $c_i= (M/p_i)^{-1}\pmod{p_i}$.}
	\end{algorithmic}
\end{algorithm}

Assume that we have constructed a pair $(x,M)$ as above. We may choose $x$ in the range $\{-\lceil M/2\rceil,\ldots,\lfloor M/2\rfloor\}$ to keep $|x|/M$ as small as possible. Then we have
\begin{equation}
	|f(x)|/M=|(x^2 + 2x\lceil \sqrt{N}\rceil + (\lceil \sqrt{N}\rceil)^2-N)|/M
	\ll  \frac{M}{4} + \lceil \sqrt{N}\rceil + \frac{2 \sqrt{N} + 1}{M} 
\label{eq:bound1}
\end{equation}
as the worst-case upper bound for our smoothness candidates. That is quite good and, assuming a reasonable size of $M$, on average quite a bit smaller than the smoothness candidates obtained in the Quadratic Sieve (see (\ref{eq:bound})). In practice, many of these candidates will actually be slightly below $\sqrt{N}$. However, the average number of different pairs $(x,M)$ we would have to test for smoothness until finding a smooth relation is still substantial, particularly for inputs with 60 digits or more. We may apply the \textbf{smooth\_batch} procedure (see \Cref{sec:rel_work}) to test a lot of candidates at once. Still, for large enough numbers, the advantage of fast sieving will outweigh the advantage of slightly smaller smoothness candidates.

Luckily, we can do better than that. Again, what we want is to efficiently construct pairs $(x,M)$ such that $M$ divides $f(x)$ and $|x|/M$ is as small as possible. Let us start with a random pair $(x,M)$ as constructed above. We now consider the values $x_j=x+jM$ for small values $j\in\Z$. It is clear that $M\mid f(x_j)$ for all $j$. What we are looking for are values $j$ for which we know that $f(x_j)$ will be divisible $\emph{not only}$ by $M$, but also by larger primes in the factor base, i.e., primes in $\lp$. In order to find such values, we compute
\begin{equation}
r_1:=(s_{i,1}-x)\cdot M^{-1} \pmod{p_i} \text{\, and \,} r_2:=(s_{i,2}-x)\cdot M^{-1} \pmod{p_i}
\label{eq:trans}
\end{equation}
for all $p_i\in\lp$. It is easy to see that $r_1$ and $r_2$ are the residues of those values $j\in\Z$ for which $p_i\mid f(x_j)$. We also compute $r_1-p_i$ and $r_2-p_i$, and store these values in a list. Having finished the computation of these lists for all primes in $\lp$, the next step is to search for \emph{collisions}, i.e., values that occur in more than one list, for different primes. Assume that $\alpha$ occurs in the lists of the primes $p$, $q$ and $r$, then $Mpqr\mid f(x_\alpha)$. We hence define a new candidate pair $(x_\alpha,Mpqr)$. Assume that $r$ is the largest of the three primes. Then $|\alpha|<r$, and we have
\begin{align*}
	|f(x_\alpha)|/(Mpqr)&=|((x+M\alpha)^2 + 2(x+M\alpha)\lceil \sqrt{N}\rceil + (\lceil \sqrt{N}\rceil)^2-N)|/(Mpqr)\\
	&<|((M(\alpha+1))^2 + 2(M(\alpha+1))\lceil \sqrt{N}\rceil + (\lceil \sqrt{N}\rceil)^2-N)|/(Mpqr)\\
	&\ll  \frac{Mr}{pq} + \frac{2\lceil \sqrt{N}\rceil}{pq} + \frac{2 \sqrt{N} + 1}{Mpqr}.
\end{align*}
Compared to (\ref{eq:bound}) and (\ref{eq:bound1}), this value is significantly smaller, improving the chance for it to factorize completely over the remaining primes in $\mathcal{S}$. In addition, it may also still be divisible by smaller primes in $\lp$ whose residues in the computed lists are not equal, but congruent to $\alpha$.
The current implementation of SSS is built for the task of quickly generating candidate pairs $(\overline{x},\overline{M})$ of the described form. $\overline{M}$ is composed of the product $M$ of several small primes in $\mathcal{S}$ and at least three large primes in $\lp$, and $\overline{x}$ is of the form $x+M\alpha$ for suitable $\alpha\in\Z$. Having collected a certain number of these pairs, we apply the \textbf{smooth\_batch} procedure with a substantially increased chance of finding smooth relations.  

The core idea of SSS is the search for suitable smoothness candidates by looking for collisions instead of sieving. In \Cref{sec:imp}, we discuss our current implementation in greater detail. There are a few things that we do differently than in the explanation above. For example, we do not really need to keep track of the large primes in $\lp$ that divide $f(\overline{x})$. It is enough to know that they exist. So we may simply store the residues $r_1, r_1-p_i,  r_2, r_2-p_i$ in a hash table or dictionary and count them directly as they are computed. All we need are the resulting values $\alpha$ and their counts, since we are only interested in those values whose count exceed a certain threshold, e.g. $3$. The information about which primes divide the resulting $f(x_\alpha)$ is lost, but we do not need it until the linear algebra step.

Another implementation detail concerns the computation of the residues $r_1$ and $r_2$, particularly the inversion of $M$ modulo the primes in $\lp$. This is a quite expensive step and similar to switching to another polynomial in MPQS or SIQS. However, we have found a remedy for this problem that will also be explained in the course of the following subsection.

\subsection{Current implementation}\label{sec:imp}
Let us now discuss our implementation\footnote{Repository: \url{https://github.com/sbaresearch/smoothsubsumsearch}} and the used parameters. For the sizes $m$ and $n$ of the factor base $\mathcal{F}$ and the small factor base $\mathcal{S}$, we aim\footnote{Our implementation considers the first $2m$ primes for $\mathcal{F}$ and then removes those $p$ with $(N|p)=-1$ from the factor bases. Since $N$ is a quadratic non-residue for about half of these primes, this leads to slightly varying, but overall similar cardinalities as those given in \Cref{tab:fb}.} at the values given in \Cref{tab:fb}. Here, $m$ is chosen similarly as for SIQS in various factorization tools. In particular, the choices for $m$ between 26 and 100 digits are taken from the source code of the \emph{primefac} package in Python, which uses similar parameters as Msieve v1.52 and is one of the implementations used in our experiments in \Cref{sec:experiments}. In order to find a good choice for $n$, which is a new parameter in our algorithm, we conducted a number of experiments. In the end, $n=m/5$ showed satisfying results.
\begin{table}[t]
	\centering
	\scriptsize
	\captionsetup{justification=centering}
	\caption{Choice of factor base sizes $m$ and $n$}\label{tab:fb}
	\begin{tabular}{lll}
		\toprule
		Digits     &     m & n\\
		\midrule
		$\hspace{10pt}\leq 18$ & 60 & 12\\
		$19 - 25$ & 150 & 30\\
		$26 - 34$ & 200 & 40\\
		$35 - 36$ & 300 & 60\\
		$37 - 38$ & 400 & 80\\
		$39 - 40$ & 500 & 100\\
		$41 - 42$ & 600 & 120\\
		$43 - 44$ & 700 & 140\\
		$45 - 48$ & 1000 & 200\\
		$49 - 52$ & 1200 & 240\\
		$53 - 56$ & 2000 & 400\\
		$57 - 60$ & 4000 & 800\\
		$61 - 66$ & 6000 & 1200\\
		$67 - 74$ & 10000 & 2000\\
		$78 - 80$ & 30000 & 6000\\
		$81 - 88$ & 50000 & 10000\\
		$89 - 94$ & 60000 & 12000\\
		$95 - 100$ & 100000 & 20000\\
		\bottomrule
	\end{tabular}
\end{table}

Let us start by discussing the precomputations that we conduct prior to our main loop for finding smooth relations.
\begin{enumerate}[label=(\roman*)]
	\item{Compute the factor base $\mathcal{F}$ and small factor base $\mathcal{S}$.}
	\item{Compute the roots of $f$ modulo $p_i$ for the odd primes $p_i$ in $\mathcal{F}$.}
	\item{Compute a product tree of $\mathcal{F}$ as preparation for \textbf{smooth\_batch}\footnote{We modified the code from \url{https://facthacks.cr.yp.to} for our implementation.} (see \Cref{sec:rel_work}).}
	\item{We want to avoid the inversion in the coefficients $c_i= (M/p_i)^{-1}\pmod{p_i}$ used in Step 2 of the function \textbf{get\_x}. Let $\mu$ be the product of all primes in $\mathcal{S}$ and let $\gamma_i:=(\mu/p_i)^{-1}\pmod{p_i}$. Then it is easy to check that, for every $M\mid \mu$ and every $i=1,\ldots,n$, we have
		\[
		\frac{\mu\gamma_i}{p_i}\equiv \frac{Mc_i}{p_i}\mod M.
		\]
		We hence precompute the global coefficients $\Lambda_i:=(\mu\gamma_i)/p_i$ for $i=1,\ldots,n$ once and may then use them in the main loop, even though the modulus $M$ is constantly changing.} 
	\item{In Step (ii), we computed the roots $\{s_{i,1}, s_{i,2}\}$ of $f$ modulo the odd primes in the small factor base $\mathcal{S}$. In addition, we also compute the differences 
		\[
		\Delta_i:=\Lambda_i\cdot(s_{i,2}-s_{i,1}).
		\]
		If we then want to switch from a solution $x\,\widehat{=}\,(x_1,\ldots,x_{i-1},1,x_{i+1},\ldots,x_n)$ to the solution $(x_1,\ldots,x_{i-1},2,x_{i+1},\ldots,x_n)$ or vice versa, we may easily do that by computing $(x+\Delta_i)\pmod M$ or $(x-\Delta_i)\pmod M$, respectively.}
\end{enumerate} 
Equipped with the quantities $\Lambda_i$ and $\Delta_i$, we may now move on to the core of the SSS algorithm, which is the search function. This function runs in a repeated loop until we have found enough smooth relations. Let $\mathcal{R}$ and $\mathcal{P}$ be the sets in which we save full and partial relations, respectively.

\begin{algorithm}
	\caption*{\textbf{search$(k,\mathcal{R},\mathcal{P})$}}
	\begin{algorithmic}[1]
		\State{Choose $k$ random indices in $\{1,\ldots,n\}$ and store them in a set $\mathcal{I}$.}
		\State{Compute the list $\mathcal{M}$ of all primes $p_i\in\mathcal{S}$ with $i\in\mathcal{I}$, and their product $M$.}
		\State{Compute the inverses $\mu_p:=M^{-1}\pmod p$ for all $p\in\lp$.}
		\State{For $i=1,\ldots,n$, let $x_i=1$ if $i\in\mathcal{I}$, and set $x_i=0$ otherwise. Compute $x:=\textbf{get\_x}(x_1,x_2,\ldots,x_n)$.}
		\For{$i\in\mathcal{I}$}
			\State{Update $x\gets(x+\Delta_i)\pmod M$ and initialize a list $\mathcal{B}$.}
			\For{$p\in\lp$}
				\State{Compute $r_k:=(s_k-x)\cdot \mu_p \pmod p$ for the roots $s_1$, $s_2$ of $f \pmod p$.}
				\State{Add $(p,r_1,r_2)$ to $\mathcal{B}$. }
			\EndFor
			\State{Initialize a set of candidates $\mathcal{C}=\emptyset$.}			
			\For{each prime $q\in\mathcal{M}$}
				\State{If $q=p_i$, \textbf{continue}}
				\State{Compute $m:=M/q$ and initialize a list $\mathcal{K}$.}
				\For{$p\in\lp$}
					\State{Compute $\alpha_k:=q\cdot r_k\pmod p$ for $k=1,2$, where $(p,r_1,r_2)\in\mathcal{B}$.} 
					\State{Add $\alpha_1,\alpha_1-p, \alpha_2, \alpha_2-p$ to $\mathcal{K}$}.
				\EndFor
				\For{each $\alpha\in\mathcal{K}$ that occurs at least three times}
					\State{Compute $\overline{x}:=x+\alpha m$ and add $f(\overline{x})/m$ to $\mathcal{C}$.}
				\EndFor
			\EndFor
			\State{Run \textbf{smooth\_batch}($\mathcal{F}$, $\mathcal{C}$). Let $g_1,\ldots, g_\ell$ be the returned values.}
			\For{$i=1,\ldots,\ell$}
				\State{If $g_i=1$, add the corresponding $\overline{x}$ to $\mathcal{R}$.}
				\State{If $1<g_i<128\cdot p_m$, add $(\overline{x},g_i)$ with the corresponding $\overline{x}$ to $\mathcal{P}$.}
			\EndFor
		\EndFor
	\end{algorithmic}
\end{algorithm}

The parameter $k$ in the search function and the size $n$ of the small factor base together have a direct impact on the average magnitude of the values $M$. Similarly to our approach for choosing $n$, we have conducted a number of experiments to derive a suitable choice for $k$. In the current implementation, we are using $k=6$. In Step 3, we then compute all inverses of $M$ modulo the primes in $\lp$. As already mentioned in the end of the last subsection, we have found a way to only have to do this once for every full loop of \textbf{search}. In the application of \textbf{get\_x} in Step 4, we use the precomputed values $\Lambda_i$. Here, we are starting with the solution $(x_1,\ldots,x_n)$ that consists of only 0's and 1's. The remainder of the algorithm is a large for-loop. In Step 6, we make sure that we are considering a different value for $x$ in each run. By adding $\Delta_i$ to the previously used $x$, we exchange the corresponding 1 in the solution vector by a 2, which avoids redundancies in the found relations. The Steps 6 to 9 compute the transformations $r_1$ and $r_2$ (see (\ref{eq:trans})) of the roots of $f\pmod p$ for $p \in \lp$. Instead of checking for collisions\footnote{In our actual implementation, we also use $r_1$ and $r_2$ from Step 8 for finding collisions. We have omitted it here because it would have made the description of \textbf{search} unnecessarily complicated.} as discussed in \Cref{sec:bi} and moving on to the next $i\in\mathcal{I}$, we save the values for later usage in the inner loop, starting at Step 11. In this loop, we are performing the search for collisions \emph{not} for our initial candidate pair $(x,M)$, but for the candidate pairs $(x,M/q)$ for all primes $q$ dividing $M$. Using the values in the precomputed list $\mathcal{B}$, Step 15 computes the respective transformation of the roots of $f\pmod p$ much more efficiently. Namely, we have
\[
(s_k-x)\cdot (M/q)^{-1} \equiv (s_k-x)\cdot \mu_p\cdot q \equiv q\cdot r_k \mod p.
\]
So the values for the list $\mathcal{K}$ that we want to check for collisions may be computed by simple multiplications of two small numbers $q$ and $r_k$ modulo $p$. This makes the method much more efficient than if we would have to compute inverses in the innermost loop.

Step 12 removes a redundancy in the computation, as the change applied to $x$ in Step 6 does not affect the values of $x\pmod{M/q}$ in the case where $q=p_i$. In Step 17, we need to find those $\alpha$ that occur more than three times in $\mathcal{K}$. There are many different ways to achieve this, from counting the values on the fly in Step 15 to hash tables and dictionaries\footnote{In Python, we may use the Counter method from the Python Standard Library collections (\url{https://docs.python.org/3/library/collections.html\#collections.Counter}).}.
Finally, we apply \textbf{smooth\_batch} to $\mathcal{C}$ and check if the resulting non-smooth parts of any of the candidates equals $1$, in which case we have found a full relation, or is smaller than $128\cdot p_m$, in which case we have found a partial relation (see \Cref{sec:rel_work}). 

The Steps 1 to 4 may be considered as a \emph{global search}, in which we end up with an initial candidate pair $(x,M)$. The large for-loop spanning from Step 5 to Step 22 may be considered as a \emph{local search}, in which we refine the initial candidate pair and obtain arguments $\overline{x}$ with a higher probability of leading to smooth values $f(\overline{x})$. We note that, after the random choice of indices in the global search in Step 1, the local search and the remainder of the algorithm are fully deterministic.
We have to repeat the function \textbf{search$(k,\mathcal{R},\mathcal{P})$} until we have found enough full and partial relations. As soon as this is the case, we proceed to the linear algebra stage of the algorithm, which is the same as for other Dixon-type algorithms. In addition to the potential optimizations that will be described in \Cref{sec:opt}, there are minor tweaks and variations of the algorithm that make it run a bit faster. Let us discuss two of them:
\begin{enumerate}
	\item{Using a different polynomial than $f(x)=(x+\lceil \sqrt{N} \rceil)^2-N$. While SSS does not really gain anything by changing polynomials during the algorithm like in SIQS, the search seems to run faster for certain polynomials compared to others. In particular, consider $f_\gamma(x)=(x+\lceil \sqrt{\gamma N} \rceil)^2-\gamma N$ for small choices of $\gamma\in\N$. A well-known approach for finding a good choice of $\gamma$ is the Knuth-Schroeppel function (\cite[p.335]{Sil}). However, we have not used this improvement in the experiments in \Cref{sec:experiments}.}
	\item{Let $\eta$ be the product of all primes in $\mathcal{F}$, which is computed in the product tree during the initialization (see Step (iii) of the precomputation steps). Before starting the search loop, we multiply $\eta$ by powers of the primes in the factor base until each prime power dividing $\eta$ is larger than $2^{15}$. This allows to reduce the exponent $e$ in Step 3 of \textbf{smooth\_batch}. In fact, our implementation does not run Step 3 at all. We might miss a few smooth relations, but the runtime savings due to the omitted repeated squarings make up for that. In slightly different form, this improvement is already mentioned in Bernstein's original paper on the smooth-batch algorithm. In fact, a long list of other ideas for speedups of the procedure can be found in \cite[Section 3]{Ber}, so we suspect that there is room for further improvement.}
\end{enumerate}

\subsection{Comparison to SIQS}\label{sec:siqs}
Let us now compare SSS to the Quadratic Sieve, and in particular to the self-initializing version SIQS that has been described in \Cref{sec:rel_work}. In the experiments in \Cref{sec:experiments}, we will see that SSS is consistently faster than SIQS by a factor between 5 and 7. While it does not appear feasible to rigorously prove a runtime improvement by a constant between such complex procedures, we will consider some arguments that provide theoretical and heuristical plausibility for the improvement by the observed factor.

Starting with the similarities, the Steps 1 to 9 in \textbf{search} are rather close to what happens in SIQS when a new polynomial $f_{a,b}$ is chosen. Our choices of $M$ and $x$ resemble the choices of $a$ and $b$, and the computation of the values $r_k$ in Step 8 corresponds to the computation of the roots of $f_{a,b}$ modulo the primes in $\lp$. While SIQS conducts these computations modulo the complete factor base $\mathcal{F}$, we restrict our attention to $\lp$. There are also differences in how we choose $M$ compared to how $a$ is typically chosen. Nevertheless, the setup and preparations for what happens later in both algorithms are somewhat similar, and it seems fair to assume that there are no major deviations in the runtimes of these parts of the techniques. The main structural difference is in the Steps 10 to 22 of \textbf{search}. In particular, we find our smoothness candidates by searching for collisions in the list $\mathcal{K}$, while SIQS uses a sieving procedure on an interval of certain length $2L$. We will now compare the two approaches for finding suitable candidates $x+\alpha M$ with $|\alpha|$ bounded by $p_n$, the largest prime in the small factor base $\mathcal{S}$. For the sake of simplicity, let us consider just the positive interval $[0,p_n]$. In SIQS, we sieve this interval in the manner described in \Cref{sec:rel_work}, which takes around
\[
2\cdot\sum_{p\in\mathcal{F}} \lceil p_n/p\rceil \geq 2\cdot(p_n\sum_{p\in\mathcal{S}} 1/p + |\mathcal{F}|-|\mathcal{S}|)
\]
operations. The factor $2$ comes from the fact that we have to do this two times, once for each root modulo each prime. In SSS, we have to compute the values $\alpha_k$ in Step 15, which are two operations (one for each root) per prime in $\lp$. Overall, this makes
\[
2\cdot (|\mathcal{F}|-|\mathcal{S}|)
\]
operations. One can now compare these two terms asymptotically and numerically. Asymptotically, the Prime Number Theorem implies that the SIQS term is around $O(p_n\log\log p_n)$ and that the SSS term is around $O(p_m/\log p_m)$, where $p_m$ is the largest prime in the factor base. So one would expect that, numerically, the difference between the two quantities grows larger. We verified this by computing the terms for a few examples. For inputs with around 30-70 digits, the SIQS term was usually about 3 to 8 times larger than the SSS term. For inputs with 75-100 digits, it was usually about 8 to 13 times larger.

Of course, this is not the complete picture. One might argue that the sieving operations in SIQS are fixed-precision additions, whereas the main operations in SSS are multiplications modulo the primes in $\lp$. While that is true, we point out that SIQS has to compute new roots of the polynomials every time it starts with a new interval, and this happens also by performing multiplications modulo the primes in the factor base. SIQS then proceeds with the sieving step. SSS, on the other hand, does nothing else but computing the multiplications. It basically omits the sieving step. 

Another objection might be that SIQS finds (almost) all candidates of the form $x+\alpha M$ with $\alpha<p_n$ that lead to smooth relations, while SSS only finds those for which the corresponding y-value has at least three prime factors in $\lp$. However, it appears that the number of relations satisfying this condition is quite large. If anything, it may make sense to further strengthen the condition, not weaken it. Depending on the set threshold in the sieving process, the resulting candidates in SIQS will have a high likelihood of leading to either full or partial relations. The condition of having at least three large prime divisors in SSS leads to more ``false positives'' that have to be tested and rejected in the smooth-batch procedure. This fact becomes more important for larger inputs with 70 digits or more (see \Cref{fig:fps} in \Cref{sec:experiments}), where the higher cost of the smooth-batch procedure in SSS partially compensates the growing amount of operations in SIQS mentioned above.

We now want to stress a few additional aspects that work to the favor of SSS.
\begin{enumerate}
	\item{Working modulo $M/q$ for the different prime factors $q$ of $M$ in the innermost loop allows us to switch to a new search domain with essentially zero cost. SSS focuses on relatively small values for $\alpha$, particularly those in the interval $[-p_n,p_n]$. SIQS has to sieve a larger interval of length $2L$ for each polynomial to reduce the cost resulting from switching between polynomials\footnote{For instance, consider inputs with 60 digits. Then $p_n$ is around 14000, while a usual choice for $L$ for inputs of this size is $196608$.}. The probability of finding smooth relations decreases sharply for growing values of $\alpha$.}
	\item{Quite regularly, SSS will find smooth relations from $\alpha$'s with absolute value larger than $p_n$. If three or more of the larger primes in $\lp$ divide a candidate value, and the corresponding residue is larger than $p_n$, this will also be detected.}
	\item{SSS is arguably simpler than SIQS. Its core functions can be implemented in around 100 lines of code, and the results on smaller numbers show that there appears to be a lower amount of computational overhead.}
\end{enumerate} 
With all this in mind, the results of our experiments (an improvement by a factor of around 10 for 30-40 digits, and by a factor of 5 to 7 for 45-100 digits) might be explained as follows: For smaller numbers (30-70 digits), the improvement is a combination of a reduced amount of overhead together with a reduced factor of around 3 to 8 of main SSS operations compared to main SIQS operations. For larger numbers, the improvement is a result of a smaller average size of $\alpha$ and a reduced factor of around 8 to 13 of main SSS operations compared to main SIQS operations. With increased input size, the smooth-batch procedure for detecting false positives becomes a bottleneck in SSS, which compensates these effects to some extent. Overall, this leads to a consistent improvement by a factor of around 5 to 7. 

As initially mentioned, this subsection is not considered to be a rigorous proof of the observed runtime reduction. There are several factors that cannot be considered in full detail here, such as particular nuances of the implementations. However, we hope to have made a convincing case for the plausibility of the improvement.

\subsection{SSSf: A filter for smoothness candidates}\label{sec:sssf}
We have said that the smooth-batch procedure gets more expensive as the input number grows and smooth relations become rare. In practical terms, this becomes noticeable on inputs $N$ with about 70 digits or more. Before presenting our experiments, we  propose an improvement that works like a filter for the candidates in $\mathcal{C}$ in Step 19 of \textbf{search}. Prior to applying \textbf{smooth\_batch} with the complete factor base $\mathcal{F}$, we apply it with a subset of its smallest primes. To be more concrete, we split up $\mathcal{F}$ into two disjoint subsets $\mathcal{F}_1$ and $\mathcal{F}_2$. $\mathcal{F}_1$ contains the smallest primes in $\mathcal{F}$ such that $\rho=|\mathcal{F}|/|\mathcal{F}_1|$ for some preset proportion $\rho$, and $\mathcal{F}_2$ contains all remaining primes. We then replace Step 19 in \textbf{search} by the following procedure.
\begin{algorithm}
	\caption*{\textbf{smooth\_filter}($\mathcal{F}_1$, $\mathcal{F}_2$, $\mathcal{C}$, $\delta$)}
	\begin{algorithmic}[1]
		\State{Initialize $\mathcal{C}_f=\emptyset$}
		\State{Run \textbf{smooth\_batch}($\mathcal{F}_1$, $\mathcal{C}$). Let $g_1,\ldots, g_\ell$ be the returned values.}
		\For{$i=1,\ldots,\ell$}
			\State{If $g_i<10^{d/2-\delta}$ for the number $d$ of decimal digits of $N$, add $g_i$ to $\mathcal{C}_f$.}
		\EndFor
		\State{Run \textbf{smooth\_batch}($\mathcal{F}_2$, $\mathcal{C}_f$) and proceed with Step 20 of \textbf{search}.}
	\end{algorithmic}
\end{algorithm}

The main goal of \textbf{smooth\_filter} is to reduce the number of candidates that have to be checked in the smooth-batch procedure for divisibility by the complete factor base. In order to achieve that, we first apply \textbf{smooth\_batch}($\mathcal{F}_1$, $\mathcal{C}$) to check for divisibility by the smallest primes. If the resulting non-smooth part of a $g_i$ is not significantly smaller than $\lceil\sqrt{N}\rceil$ in Step 4, we directly discard this candidate.

We denote SSS applied with this modification by SSSf. While SSSf is in fact slower than SSS for inputs $N$ with up to around 70 digits, it shows an improved performance on inputs with 75-100 digits. Besides the already discussed parameters $m$, $n$ and $k$, SSSf comes with two additional parameters $\rho$ and $\delta$. In our current implementation, we are using $\rho=10$ and $\delta=5$. In addition, we are using $k=7$ instead of the choice of $k=6$ we have used in SSS. However, these values may not be optimal, and it certainly makes sense to conduct more experiments for fine-tuning them with regards to the size of the input number.

\section{Experiments}\label{sec:experiments}
The approach described in \Cref{sec:sss} has been implemented in the programming language Python, a freely available high-level language that emphasizes code readability. The reader may therefore easily reproduce all experimental results presented in this section. On the other hand, Python is considerably slower than other programming languages (e.g. C), which restricted the size of the input numbers $N$ we could factorize in our experiments.
In the range from 30 to 70 decimal digits, it was possible to perform full factorization tests on random semiprimes. For larger inputs ranging from 75 to 100 digits, this would have taken too much time. For such numbers, we stopped the algorithms after one hour and compared the amount of smooth relations that have been found during this time. Due to this limitation, there is a factor of randomness involved, particularly for the results on inputs with 90 to 100 digits. However, it allows us to observe a general tendency of the runtime of the algorithms and their development for growing input numbers.

\begin{figure}[ht]
	\centering
	\begin{subfigure}{.48\textwidth}
		\centering
		\includegraphics[width=1.07\linewidth]{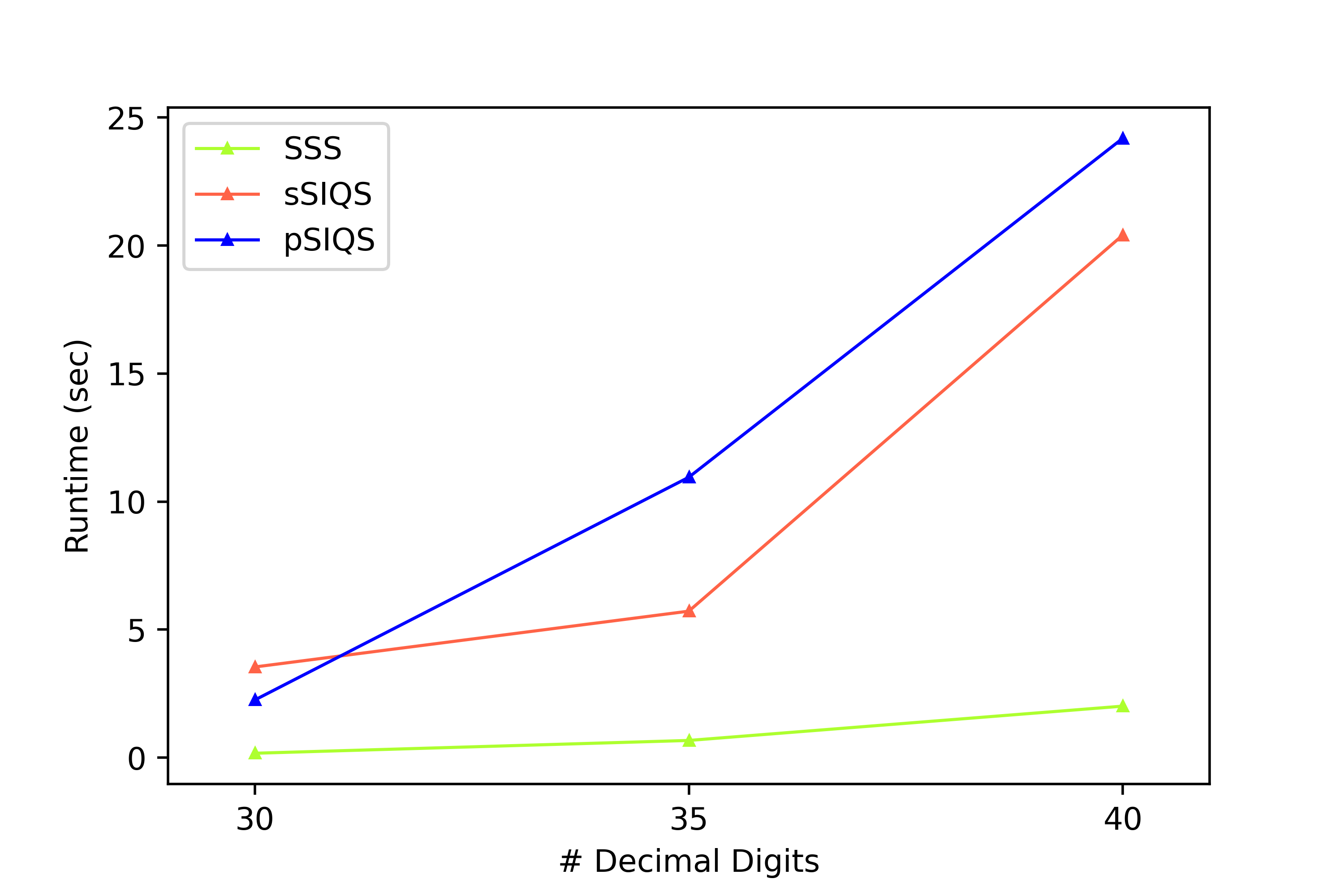}
		\caption{30 to 40 digits} \label{fig:30_40}
	\end{subfigure}\quad
	\begin{subfigure}{.48\textwidth}
		\centering
		\includegraphics[width=1.07\linewidth]{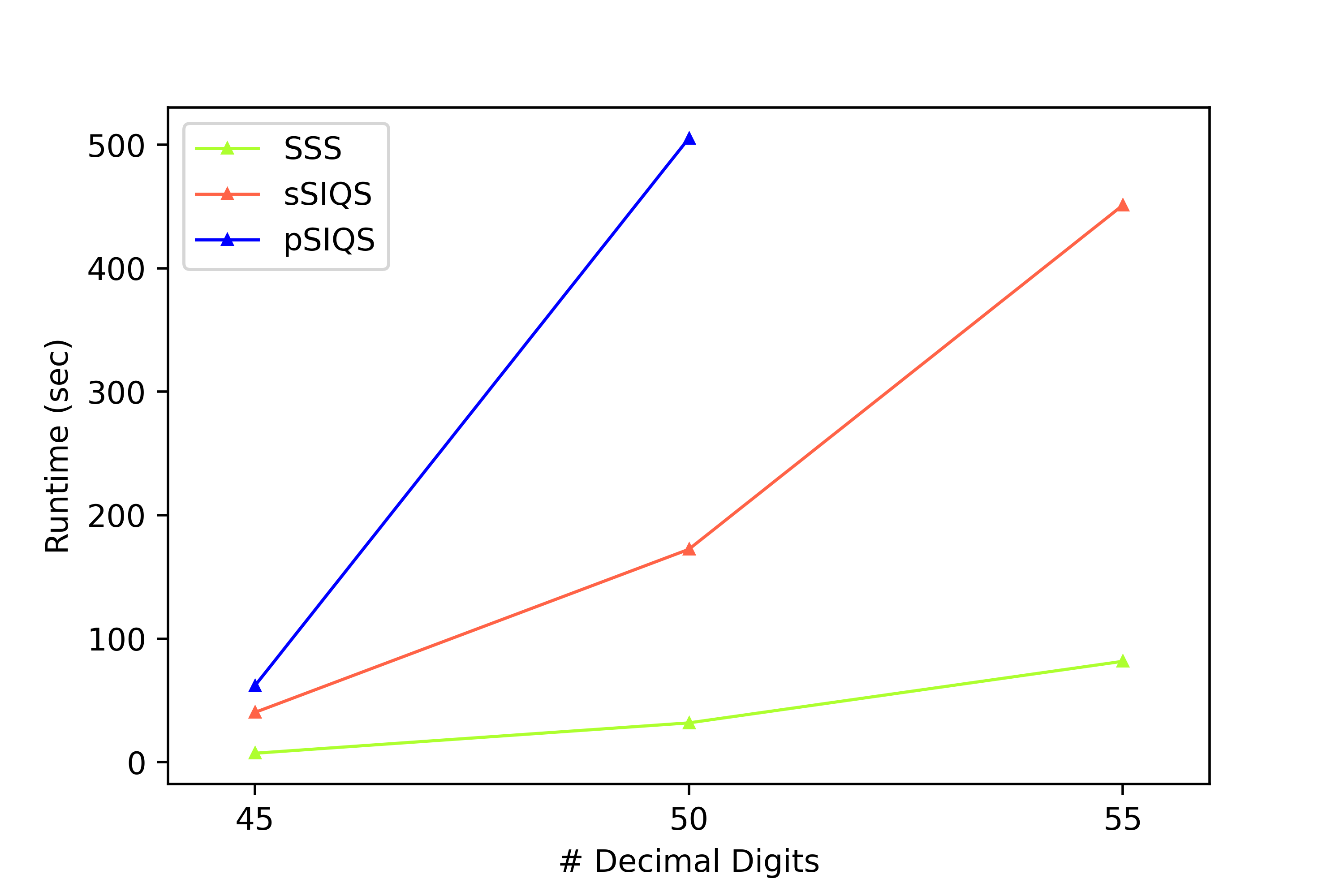}
		\caption{45 to 55 digits} \label{fig:45_55}
	\end{subfigure}
	\begin{subfigure}{.48\textwidth}
		\centering
		\includegraphics[width=1.07\linewidth]{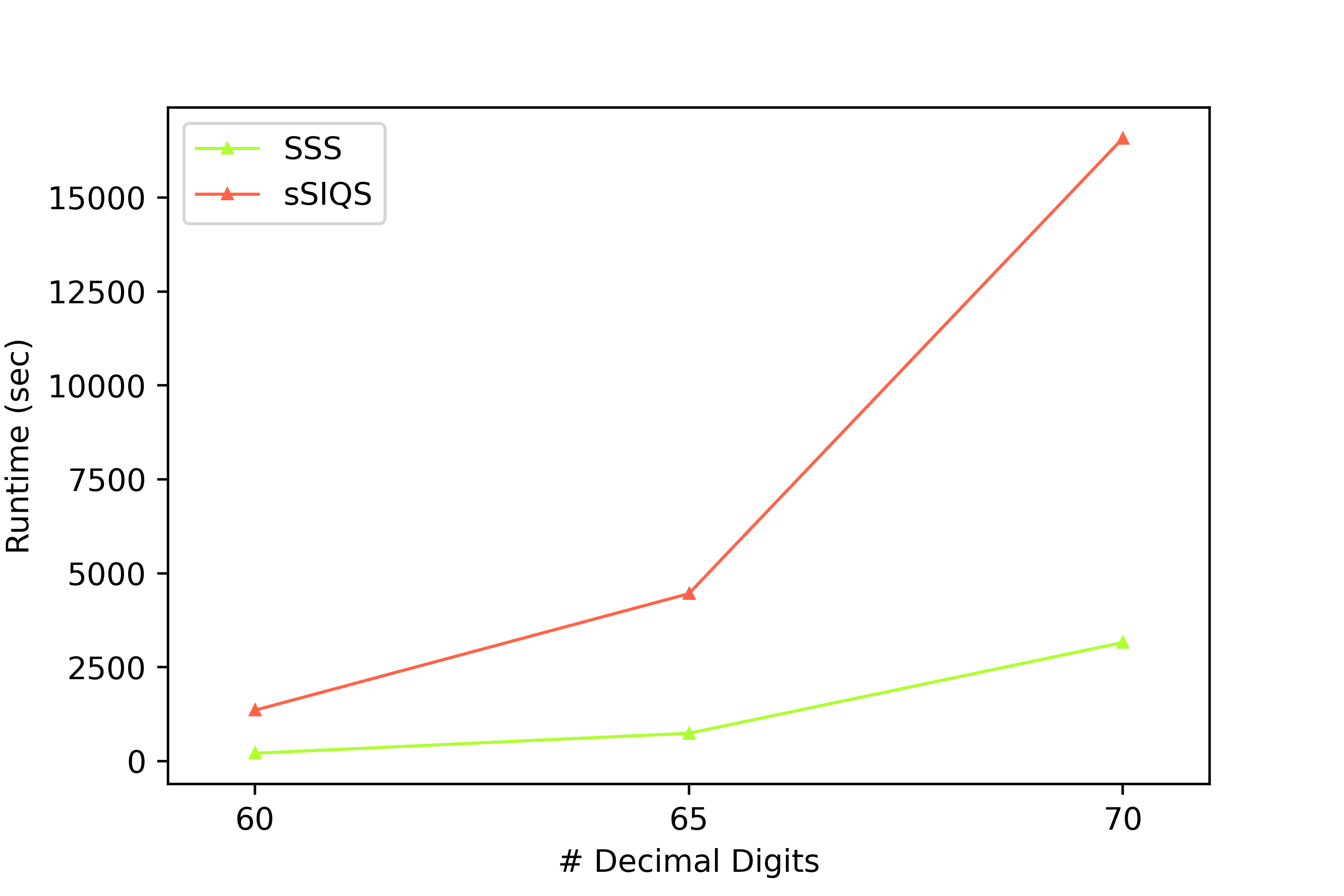}
		\caption{60 to 70 digits}
		\label{fig:60_70}
	\end{subfigure} 
\caption{Mean runtime complexity results as shown in \Cref{tab:runtime}}
\label{fig:runtime}
\end{figure}

In order to compare SSS to SIQS, we have used currently available implementations of SIQS in Python libraries. The first we found is a recent implementation in \emph{sympy}\footnote{\url{https://docs.sympy.org/latest/modules/ntheory.html}}, a well-known and widely used library for symbolic mathematics. The second implementation we found is part of \emph{primefac}\footnote{\url{https://pypi.org/project/primefac/}}, a package focused on primality tests and integer factorization. It is noted in the documentation that their SIQS implementation is taken mostly verbatim from another package, \emph{PyFactorise}\footnote{\url{https://github.com/skollmann/PyFactorise}}. While we have found other implementations of the (Self-initializing) Quadratic Sieve in Python, these two were most refined and competitive in our experiments. Throughout the section, we will denote sympy's implementation by sSIQS and the implementation by primefac/PyFactorise by pSIQS.
With regards to the comparability of the approaches, we point out the following.
\begin{itemize}
	\item{We have applied all three implementations with the same code for the second phase (the linear algebra stage\footnote{In general, this stage contributes only a tiny fraction to the complete runtime.}) to make sure that any runtime differences only result from the performance in the first phase, where smooth relations are collected. We utilized a Gaussian elimination algorithm over $\mathbb{F}_2$ as it is implemented in primefac. The approach is based on \cite{Koc}.}
	\item{For all approaches, we used a factor base size $m$ of $\mathcal{F}$ close to the values in \Cref{tab:fb}. For the length $L$ of the sieved intervals in the SIQS implementations, we used the standard settings given in the source code of pSIQS.}
	\item{Both sSIQS and pSIQS do not use the Knuth-Schroeppel function to find a suitable multiplier $\gamma$ for $N$. As already discussed at the end of \Cref{sec:imp}, we also did not make use of this improvement.}
	\item{In pSIQS, the large prime variant (see end of \Cref{sec:rel_work}) is not implemented. This has to be kept in mind when considering the results, and is also the reason why we did not run pSIQS on inputs with 55 digits or more.}
	\item{We used SSS with the bound $128\cdot p_m$ for saving partial relations (see Step 22 of \textbf{search}). The same bound is used in sSIQS.}
\end{itemize}
Let us now discuss our experimental setup. We consider input numbers $N$ with $d=30, 35, 40, \ldots, 95, 100$ digits. For each of these input sizes, we did the following.
\begin{enumerate}
	\item{Generate a random semiprime $N$ with $d$ digits, such that its two prime factors are about the same size.}
	\item{Apply sSIQS to factorize $N$. If $d\leq 50$, apply pSIQS to factorize $N$.}
	\item{If $d\leq 70$, apply SSS to $N$. For $d=75,80,85,90,95,100$, apply SSSf to $N$.}
\end{enumerate}
Let us start with the results\footnote{All computations have been conducted on a standard laptop with AMD Ryzen 7 5800H processor (8-core 3.2 GHz) and 16GB RAM.} for $30\leq d\leq 70$. Due to the increase in general time complexity, we consider $20$ different inputs $N$ for $d\in\{30,35,40\}$, $10$ different $N$ for $d\in\{45,50,55\}$, five different $N$ for $d=60$  and three different $N$ for $d=65,70$. After we have finished all runs for a certain $d$, we computed the mean and the standard deviation of the runtimes of all methods. The results can be found in \Cref{tab:runtime}, and the mean runtimes are visualized in \Cref{fig:runtime} in form of line plots.
\begin{table}[ht]
	\centering
	\small
	\setlength\tabcolsep{4pt}
	\captionsetup{justification=centering}
	\caption{Runtime complexity in seconds: Mean $\pm$ STD\\
		Last column displays improvement factor of SSS over the fastest SIQS method}\label{tab:runtime}
	\begin{tabular}{lllll}
		\toprule
		\#Digits     &     SSS &         sSIQS &       pSIQS  & \\
		\midrule
		30 & 0.17$\pm$0.06 & 3.54$\pm$2.51 & 2.25$\pm$1.33& 13.2\\
		35 & 0.67$\pm$0.27 & 5.72$\pm$2.24 & 10.95$\pm$6.98& 8.5\\
		40 & 2.01$\pm$0.64 & 20.41$\pm$5.72 & 24.19$\pm$12.26& 10.2\\
		\midrule
		45 & 7.16$\pm$2.85 & 40.20$\pm$17.21 & 61.71$\pm$40.66& 5.6\\
		50 & 31.76$\pm$13.04 & 172.35$\pm$55.96 & 505.34$\pm$282.97& 5.4\\
		55 & 81.64$\pm$32.94 & 451.13$\pm$162.42 & -- & 5.5\\
		\midrule
		60 & 208.35$\pm$32.67 & 1352.57$\pm$176.28 & -- & 6.5\\
		65 & 740.12$\pm$120.89 & 4456.17$\pm$984.49 & --& 6.0\\
		70 & 3157.13$\pm$252.08 & 16578.18$\pm$1532.10 & --& 5.3\\
		\bottomrule
	\end{tabular}
\end{table}

For inputs with $30-40$ digits, SSS is between 8.5 and 13.2 times faster than the SIQS approaches. For inputs with $45-70$ digits, the improvement factor is in the range between 5.3 and 6.5. sSIQS and pSIQS behave somewhat similar for integers with 30 to 45 digits. For inputs with 50 digits and higher, it becomes increasingly noticeable that pSIQS does not make use of the large prime variant. Hence, we did not consider it in our experiments on numbers with 55 to 70 digits.

Let us now discuss the results for inputs with 75 to 100 digits. For each $d$ in this range, we conducted experiments on two different numbers $N$. For each $N$, we applied SSSf (see \Cref{sec:sssf}) and sSIQS and let them run for one hour each. Finally, we compared the number of found relations during this time. The results are in \Cref{tab:runtime_high}.
\begin{table}[ht]
	\centering
	\small
	\setlength\tabcolsep{4pt}
	\captionsetup{justification=centering}
	\caption{Found relations after 3600 seconds: Mean $\pm$ STD\\
		Last column displays improvement factor of SSSf over sSIQS}\label{tab:runtime_high}
	\begin{tabular}{llll}
		\toprule
		\#Digits     &     SSSf &         sSIQS  & \\
		\midrule
		75 & 12737.5$\pm$2677.5 & 1880.0$\pm$363.0 &  6.8\\
		80 & 5544.5$\pm$1102.5 & 920.5$\pm$146.5 &  6.0\\
		85 & 1790.5$\pm$539.5 & 350.0$\pm$100.0 & 5.1\\
		\midrule
		90 & 538.5$\pm$36.5 & 86.0$\pm$16.0 &  6.3\\
		95 & 356.0$\pm$4.0 & 64.5$\pm$5.5 &  5.5\\
		100 & 139.0$\pm$51.0 & 25.0$\pm$5.0 &  5.6\\
		\bottomrule
	\end{tabular}
\end{table}
 
We can see that SSSf finds between 5.1 and 6.8 times as many smooth relations as sSIQS in the same time. As already mentioned, these are by no means large scale experiments, and there is a certain amount of variance involved in the specific numbers of the table. Looking at general tendencies, however, we point out that the factors in the last column are astonishingly similar to the results we got for numbers with 45 to 70 digits in \Cref{tab:runtime}. In particular, there appears to be no substantial decrease of the advantage of SSSf over sSIQS for larger inputs. 

With regards to the setup of this experiment, one should mention that the number of found relations in SSS, SSSf and SIQS is not uniformly distributed over time. In the later stages of the process, the methods are more likely to connect two partial relations to a full relation and thereby, on average, find more smooth numbers compared to the beginning of the process. However, this is something that affects both SSSf and sSIQS in the same manner. As a consequence, it seems reasonable to conclude from these results that SSSf would finish the factorization of numbers with 75-100 digits around 5 to 7 times faster than sSIQS.
\begin{figure}[ht]
	\centering
		\includegraphics[width=.7\linewidth]{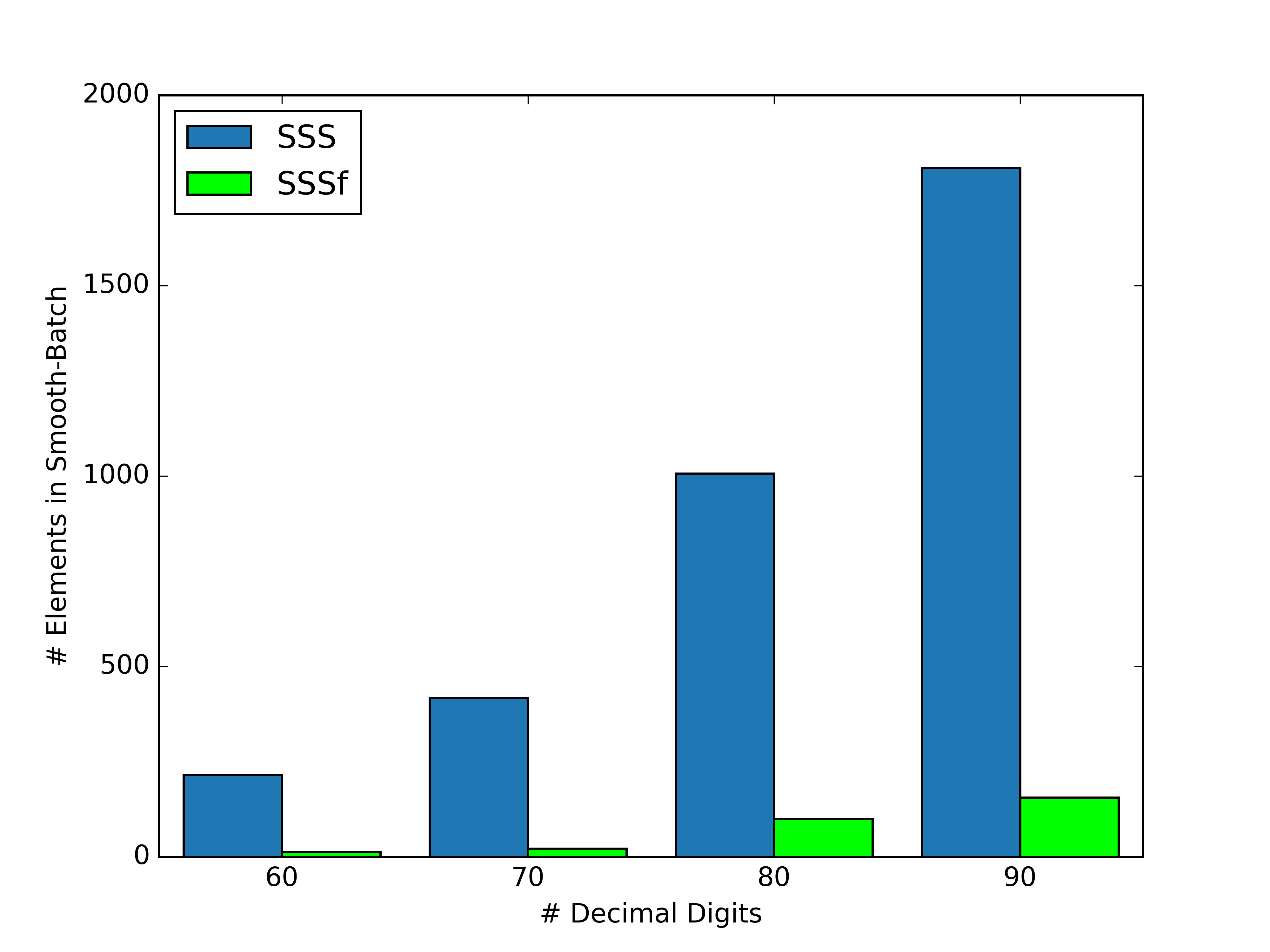}
	\caption{Number of average elements in $\mathcal{C}$ in Step 19 of \textbf{search}. We compare the size of $\mathcal{C}$ in SSS with the size of the filtered set $\mathcal{C}_f$ of candidates in in \textbf{smooth\_filter} of SSSf.}
	\label{fig:fps}
\end{figure}

\vspace{12pt}
Finally, let us consider the cardinality of the candidate set $\mathcal{C}$ to which we apply \textbf{smooth\_batch} in Step 19 of \textbf{search}. In \Cref{sec:siqs}, we have already talked about the number of false positives in SSS. While this number is rather small for inputs with 55 digits or less, it starts to increase for inputs with more than 65-70 digits. In \Cref{fig:fps}, a histogram shows the average size of $\mathcal{C}$ for inputs with 60, 70, 80 and 90 digits. We can see that the cardinality of $\mathcal{C}$ is between 1000 and 2000 for inputs with 80 to 90 digits. This is still a very small amount compared to the length of the search domain, but it constitutes an additional amount of work in \textbf{smooth\_batch} that, in this form, is not necessary in SIQS. There, the number of elements that survive the sieving is comparably small, and so is the cost for processing them, i.e., checking whether they lead to full or partial relations.

One solution for mitigating the additional cost is to change our condition for including elements in $\mathcal{C}$. For example, we could increase the number of large primes in $\lp$ that have to divide our candidates from $3$ to $4$. Indeed, this tends to drastically reduce the size of $\mathcal{C}$ to a two-digit number. However, it also reduces the number of found relations and appears to lead to an overall slower algorithm. A much better solution is the application of SSSf as presented in \Cref{sec:sssf}. We first apply the smooth-batch procedure to all candidates in $\mathcal{C}$, but only modulo a fraction (e.g. $1/10$) of the factor base. This is far less expensive and allows us to discard a lot of less-promising candidates that are not small enough afterwards. \Cref{fig:fps} shows that this filter is very effective in reducing the number of elements to which we have to apply the smooth-batch procedure modulo the complete factor base. While these steps are still more work compared to what we have to do in SIQS, it only appears to equalize the growing factor of increased main operations (see \Cref{sec:siqs}) in SIQS compared to the main operations in SSS(f) for inputs with 75-100 digits, and we still achieve an improvement by a factor of 5 to 7.

\section{Other Ideas for SSS}\label{sec:improvements}
Our chances of finding smooth numbers with SSS depend on the ratio of the solution $x$ and considered modulus $M$. The size of $x$ directly affects the size of the considered candidate $f(x)$, while $M$ divides $f(x)$ by definition and is a factor that is always smooth for all considered candidates. As a result, we would like to make $x$ as small and $M$ as large as possible. For a moment, let us forget about what we are doing in our current implementation, and consider the following three possible approaches for solving this problem:
\begin{enumerate}
	\item{For a fixed modulus $M$, try to find a related representation $(x_1,\ldots,x_n)$ such that the value $x$ returned by \textbf{get\_x}$(x_1,x_2,\ldots,x_n)$ is minimized.}
	\item{For a fixed value $x$, try to find divisors $m$ dividing $M$ such that the ratio $(x\pmod m)/m$ is minimized.}
	\item{Use a different search technique to optimize both parameters at once by minimizing the objective function 
		\[
		\psi=|x|/M.
		\]}
\end{enumerate}
\Cref{sec:lr}, \Cref{sec:sr} and \Cref{sec:ga} refer to these approaches respectively. Besides some proof-of-concept experiments, the described ideas have not been implemented yet. Testing their applicability to possibly improve SSS is a subject for future research. \Cref{sec:opt} discusses some minor speedups of our current code.

\subsection{Lattice Reduction} \label{sec:lr}
We start with the approach of fixing $M$ and finding a suitable subsum, i.e., a corresponding representation $(x_1,\ldots,x_n)$ such that $x$ returned by $\textbf{get\_x}(x_1,x_2,\ldots,x_n)$ is minimized. Let $M=\rho_1\ldots\rho_k$, where the $\rho_i$ are primes in our factor base or powers of such primes. For $\gamma_i:=(M/\rho_i)^{-1}\pmod {\rho_i}$ and $M_i:=M\gamma_i/\rho_i$, we may write the possible values $x$ returned by \textbf{get\_x} as
\begin{align}\label{eq:crem}
	s_{1,j_1}M_1+s_{2,j_2}M_2+\ldots+s_{k,j_k}M_k\pmod M,
\end{align}
where $s_{i,j_i}$ are solutions of $f(x)\equiv 0 \pmod{\rho_i}$. Our ultimate goal is to find the choice of $(s_{1,j_1},\ldots,s_{k,j_k})$ for which the value of (\ref{eq:crem}) is minimized. 

In principle, this problem is a modular version of the \emph{multiple choice subset-sum problem} (MCSS), which is discussed in \cite[Sec. 11.10.1]{KelPfePis}. We consider $k$ classes $N_i:=\{s_{i,j_i}M_i\pmod M: f(s_{i,j_i})\equiv 0 \pmod{\rho_i}\}$, and the elements of the classes are called \emph{weights}. Denoting our weights as $\omega_{i,j_i} := s_{i,j_i}M_i\pmod M$, our goal is to minimize the sum $\omega_{1,j_1}+\cdots+\omega_{k,j_k}\pmod M$. While MCSS is NP-complete, there are multiple ways to approach this problem and solve it efficiently for certain sizes of $k$ or cardinalities of the classes $N_i$. For example, a time-space tradeoff is discussed in \cite[Section 6]{Hit}, where a similar MCSS occurs in a more theoretical approach to integer factorization. There, the solution of the MCSS corresponds to the sum $S:=p+q$ of the two prime factors of semiprime numbers $N=pq$. While it would be enough to solve one MCSS instance in order to find $S$ (and thereby factor $N$), the problem of the approach in \cite{Hit} in terms of efficiency is that both $k$ and the classes $N_i$ become quite large (around the size of $\log N$). In our present setting, however, we are basically free to choose $k$ and may also control the cardinalities of the $N_i$ by deciding in advance which prime powers to include as factors of $M$. We may hence choose from the large variety of available techniques for minimizing (\ref{eq:crem}), one of which is based on lattice reduction.

It is well known that the LLL-algorithm (\cite{len0}) for lattice reduction can be used to solve a certain class of knapsack and subset-sum problems in polynomial-time (\cite{SchEuc}). In particular, consider the standard subset-sum problem
\begin{align}\label{eq:ssum}
	a_1x_1+\cdots +a_\ell x_\ell=s,
\end{align}
where the $a_i$ and $s$ are integers. The task is to find $(x_1,\ldots,x_\ell)\in\{0,1\}^{\ell}$ such that $(\ref{eq:ssum})$ is satisfied. An important quantity is the density 
$
d:=\ell/\log_2(\max_i a_i).
$
In \cite{CosJouMac} it is shown that an oracle for finding the shortest vector in a special lattice can be used to solve almost all subset-sum problems with $d<0.9408$. In practice, this oracle is replaced by an application of the LLL-algorithm. The procedure works as follows: For $n:=\lceil \frac{1}{2}\sqrt{\ell}\rceil$, we define the lattice $L\subseteq \Z^{\ell+1}$ generated by the rows $r_1,\ldots,r_{\ell+1}$ of the matrix
\[
\begin{pmatrix}
	1 & 0 & \cdots & 0 & na_1 \\
	0 & 1 & \cdots & 0 & na_2 \\
	\vdots  & \vdots  & \ddots & \vdots & \vdots  \\
	0 & 0 & \cdots & 1 & na_n \\ 
	\frac{1}{2} & \frac{1}{2} & \cdots & \frac{1}{2} & ns
\end{pmatrix}.
\]
If $(x_1,\ldots,x_\ell)$ is a solution to (\ref{eq:ssum}), the vector
$
v=\sum_{i=1}^\ell x_ir_i-r_{\ell+1}=(y_1,\ldots,y_\ell,0)
$
is an element of $L$, where $y_i\in\{-\frac{1}{2},\frac{1}{2}\}$ and, hence, $\norm{v}_2\leq \frac{1}{2}\sqrt{\ell}$. Knowing $v$, it is easy to retrieve a solution to the subset-sum problem. Therefore, the hope of applying the LLL-algorithm is that $v$ will occur in the reduced basis of $L$.

If we apply the approach of lattice reduction for standard subset-sum directly to MCSS, the multiple-choice restriction (choosing exactly one element from each class $N_i$) is not accounted for in the definition of the lattice $L$, hence too many vectors representing invalid solutions will be in the reduced basis. We now discuss an idea for solving this problem. As mentioned, we consider the classes $N_1,N_2,\ldots,N_k$, and each class contains the weights $\omega_{i,1},\ldots,\omega_{i,\kappa_i}$, where $i=1,\ldots,k$ and $\kappa_i:=|N_i|$. Moreover, let $\alpha\in\N$ be a natural number which will be specified later. Our idea is to introduce a set of dummy subset-sum problems in the last few columns of the matrix by which the lattice is defined. The purpose is to force the lattice reduction algorithms into favoring vectors that correspond to choosing exactly one weight from each class. For example, assume that $\kappa_i=2$ for every $i$ and that $s$ is the target sum of our MCSS. Moreover, let $n:=\lceil\frac{1}{2}\sqrt{k}\rceil$. In this case, the matrix defining the lattice is of the shape
\begin{equation}\label{eq:latmat}
	\begin{pmatrix}
		1 & 0 & 0 & 0 & \cdots & 0 & 0 & \alpha & 0 & \cdots & 0 & n\omega_{1,1} \\
		0 & 1 & 0 & 0 & \cdots & 0 & 0 & \alpha & 0 & \cdots & 0 & n\omega_{1,2} \\
		0 & 0 & 1 & 0 & \cdots & 0 & 0 & 0 & \alpha & \cdots & 0 & n\omega_{2,1}\\
		0 & 0 & 0 & 1 & \cdots & 0 & 0 & 0 & \alpha & \cdots & 0 & n\omega_{2,2} \\
		\vdots  & \vdots & \vdots & \vdots & \ddots & \vdots & \vdots  & \vdots & \vdots & \ddots & \vdots & \vdots \\
		0 & 0 & 0 & 0 & \cdots & 1 & 0 & 0 & 0 & \cdots & \alpha & n\omega_{k,1} \\
		0 & 0 & 0 & 0 & \cdots & 0 & 1 & 0 & 0 & \cdots & \alpha & n\omega_{k,2}\\
		\frac{1}{2} & \frac{1}{2} & \frac{1}{2} & \frac{1}{2} & \cdots 
		& \frac{1}{2} & \frac{1}{2} & \alpha & \alpha & \cdots & \alpha & ns\\
	\end{pmatrix}.
\end{equation}
So in addition to what we have seen in the matrix for the standard subset-sum problem, we add one additional column for each class, namely 
\[
(\alpha,\alpha,0,0,\ldots,0,\alpha)^T,\,\,\, (0,0,\alpha,\alpha,\ldots,0,\alpha)^T
\]
and so on. It is clear how the definition of this matrix may be generalized to classes with more than two weights in them. Let us consider this general case, and set $\ell=\sum_i \kappa_i$. If $(x_1,\ldots,x_\ell)$ is a binary vector encoding a solution to the MCSS (in the same way as described above for standard subset-sum), then the vector   
\[
v=\sum_{i=1}^\ell x_ir_i-r_{\ell+1}=(y_1,\ldots,y_\ell,0,\cdots, 0, 0)
\]
with $y_i\in\{-\frac{1}{2},\frac{1}{2}\}$ will be in the lattice $L\subseteq \Z^{\ell+k+1}$, and satisfies $\norm{v}_2=\frac{1}{2}\sqrt{k}$.

Besides searching for a short vector in the lattice $L$, we may also reduce MCSS to a specific Closest Vector Problem (CVP). This can be achieved by considering the lattice $L'$ defined by the matrix in (\ref{eq:latmat}) \emph{without} the last row, and then defining the target vector
$
t:=(0,\ldots, 0, \alpha, \alpha, \ldots  \alpha, Ns)^T.
$
We already know that the vector in $L'$ corresponding to the solution $(x_1,\ldots,x_\ell)$ is very close to $t$. In the context of our research on the idea in \cite{Hit}, we conducted proof-of-concept experiments and noticed that this CVP approach works much better in practice than the SVP approach discussed above. One apparent reason is that there can be vectors in $L$ that are still much shorter than $\frac{1}{2}\sqrt{k}$ in the euclidean norm. Most of these vectors relate to some combination of the weights in the last column of (\ref{eq:latmat}) that sums to $0$ and omit the last row. The formulation of the CVP prevents this behavior.

We may transform our modular MCSS into an ordinary MCSS by defining the class $N_{0}:=\{0,M,2M,\ldots (k-1)M\}$ and setting the approximative target sum $s=(k-1)M$. In addition, note that the LLL approach described above only works for an exact target sum. However, \cite[Section 6]{HowJou} shows how one may reduce any approximate MCSS instance to an exact one by dividing all weights by the approximate bound. Of course, there are several other aspects and details to consider, and we will elaborate on this approach in future work. One goal would be to prove an analogue of the standard subset-sum density result in \cite{CosJouMac} for our adaptation to MCSS. Such a result would also inform the choice of $M$ in our SSS setting. As already mentioned, we want to take $M$ as large as possible. Using lattice reduction for subsequently finding a relatively small value $x$ in the related subsum has the potential of increasing our chances of finding smooth values of polynomials $f$.

\subsection{Small modular roots of linear polynomials}\label{sec:sr}
Let us now assume that we fix the value $x$ (with regards to a certain modulus $M$) and try to find divisors $m$ of $M$ such that the ratio $(x\pmod m)/m$ is minimized. One technique that might allow us to use larger values of $M$ is Coppersmith's method (\cite{Cop}). In particular, consider the following result (\cite[Theorem 3]{may}).
\begin{theorem}\label{thm:copper}
Let $M$ be an integer which has an unknown divisor $m\geq M^\beta$, where $0<\beta\leq 1$. Let $0<\varepsilon< \frac{1}{7}\beta$. Furthermore, let $g(X)$ be a univariate monic polynomial of degree $\delta$. Then we can find all solutions $x_0$ for $g(X)\equiv 0 \mod m$ with 
\[
|x_0|\leq \frac{1}{2} M ^{\frac{\beta^2}{\delta}-\varepsilon}.
\]
The running time is polynomial in $\varepsilon^{-1},\delta$ and $\log M$.
\end{theorem}
Assume that there is some divisor $m\geq M^\beta$ of $M$ such that the residue $x_0$ of $x$ modulo $m$ satisfies $|x_0|\leq \frac{1}{2} M ^{\beta^2-\varepsilon}$. Then we can find $x_0$ by running \Cref{thm:copper} with the polynomial $g(X)=X-x$. Having found $x_0$, it is easy to determine the corresponding divisor $m$ and, thus, a suitable pair $(x_0,m)$ with low ratio.
Solving linear equations modulo unknown divisors has also been studied in \cite{her08} and \cite{lu15}. While this idea appears interesting in theory, it is not yet clear whether it is actually applicable to our problem. The main question is how to choose $M$ and $\beta$ to ensure the existence of pairs $(x_0,m)$ that satisfy the bounds. 

\subsection{A genetic algorithm}\label{sec:ga}
Finally, we assume an unrestricted search approach that chooses $x$ and $M$ together in an attempt to minimize the objective function $\psi=|x|/M$. For example, we could comb through different choices of $x$ and $M$, save all pairs $(x,M)$ for which $\psi$ is below a certain bound in a set of candidates $\mathcal{C}$, and regularly apply the smooth-batch procedure to $\mathcal{C}$ to find smooth values of $f(x)$. This idea was our starting point for the SSS implementation before applying the more refined collision-based approach. In any case, there may be other, more elaborated approaches to quickly generate (and identify) suitable pairs $(x,M)$. One such possibility is a \emph{genetic algorithm} (\cite{gold}). 
In general, genetic algorithms solve an optimization problem by representing the search space by a number of individuals (the ``population'') that represent possible solutions. Each individual is represented by a so-called ``chromosome'' that reflects its properties. In our case, we may just use the representation vectors $(x_1,\ldots,x_n)$ as discussed in the beginning of \Cref{sec:sss}. 
The next step is to find the best individuals in the population with regards to a fitness function. Here, we could use a version of our objective function $\psi=|x|/M$. Finally, the best $n$ individuals are selected as ``parents'' for the next generation, which is obtained by applying certain changes (the ``crossovers'' and ``mutations'') to the chromosomes. For example, we could change a random value in a chromosome (mutation), or we could take the chromosomes from two parents (say, $(x_1,\ldots, x_n)$ and $(y_1,\ldots, y_n)$) and generate a child chromosome as $(x_1,\ldots,x_{i-1},y_i,\ldots, y_n)$ for some random index $i$ (crossover). As soon as we have obtained enough chromosomes for the next generation, the whole process starts again, i.e., we choose the best individuals and apply mutation and crossover.

Using the chromosomes and the fitness function as discussed above, we could immediately apply a genetic algorithm to solve the SSS problem. However, the main difficulty in improving the current implementation is in finding suitable mutation and crossover operations. Ideally, they would incorporate our current strategy and have a tendency of generating individuals with a smaller ratio of $x$ and $M$ than their parents. It may also make sense to combine the genetic algorithm with the ideas discussed in the previous subsections.

\subsection{Further speedups}\label{sec:opt}
In addition to the Knuth-Schroeppel multiplier and Bernstein's optimizations for the smooth-batch procedure discussed at the end of \Cref{sec:imp}, there are other possible ways for improving SSS. Here are four examples.
\begin{itemize}
	\item{There are certain details of the algorithm that could be changed and played with. For example, we could use more than the $|\mathcal{I}|$ different values for $x$ in Step 6 of \textbf{search}, prolonging the time before having to initialize a new run of the function. We could also use other divisors of $M$ as elements in $\mathcal{M}$ in the for-loop in Step 11, not only its prime factors.}
	\item{We may control the batch size (i.e., the cardinality of $\mathcal{C}$) in \textbf{smooth\_batch}. At the moment, we just apply the smooth-batch procedure at the end of each run of the main loop in \textbf{search}. It might be better to collect more candidates until we have reached a certain number in $\mathcal{C}$ or $\mathcal{C}_f$.}
	\item{One can certainly fine-tune the currently used parameters. For example, our current choices of $\rho$ and $\delta$ in SSSf are most likely not optimal for larger numbers. It may also make sense to try other values for $m$ and $n$.}
	\item{Due to the simplicity of SSS, the whole procedure is cut out for parallelization. Different runs of \textbf{search} are independent from one another, so we may easily use more than one processor for the collection of smooth relations. In addition, the main loop in \textbf{search} is also parallelizable.}
\end{itemize} 

\section{Summary}\label{sec:conclusion}
This paper presented Smooth Subsum Search, a new heuristic search approach for finding smooth values of polynomials. We have applied SSS as part of an integer factorization algorithm. For inputs between 30 and 100 digits, we compared SSS to the Self-initializing Quadratic Sieve. Our results show that SSS runs around 10 times faster for inputs with 30 to 40 digits, and 5 to 7 times faster for inputs with 45 to 70 digits. SSSf, an improved version of the original procedure for the application on larger numbers, is also 5 to 7 times faster in finding relations for inputs with 75 to 100 digits. In addition, we have presented three ideas for possible variants and further speedups in \Cref{sec:improvements}. Our future research will concern a detailed investigation of these approaches. In particular, the lattice-based approach in \Cref{sec:lr} for solving MCSS problems could have other applications besides integer factorization. We will also work on the optimizations of the current implementation discussed in \Cref{sec:opt}. Finally, we intend to explore further applications of SSS in other algorithms that also depend on finding smooth values of polynomials, such as the Number Field Sieve.

\subsection*{Implementation} The reader may find our implementations of SSS and SSSf, the code for the runtime experiments conducted in \Cref{sec:experiments}, the inputs $N$ used in the experiments, as well as the individual results for each $N$ here:  \url{https://github.com/sbaresearch/smoothsubsumsearch}.


\end{document}